\documentclass[12pt,a4paper]{article}

\usepackage{pdflscape}
\usepackage{amsmath}
\usepackage{amssymb}
\usepackage{latexsym}
\usepackage{srcltx}
\usepackage{graphics}
\usepackage{color}
\usepackage{epsfig}
\usepackage{color}
\usepackage{graphicx,epstopdf}

\newcommand{\re}{{\mathbb R}}

\newcommand{\z}{{\mathbb Z}}

\newcommand{\bu}{{\boldsymbol{u}}}

\newcommand{\Tor}{\mathbb{T}}

\newtheorem{theorem}{Theorem}
\newtheorem{prop}{Proposition}

\newtheorem{remark}{Remark}

\date{}

\begin{document}

\author{Vladimir Yu. Protasov 
\thanks{University of L'Aquila (Italy);  {e-mail: \tt\small
vladimir.protasov@univaq.it}
}}

\title{Closed simple geodesics on a polyhedron 
 }

\maketitle

\begin{abstract}

It is well-known that every isosceles  tetrahedron (disphenoid) admits 
infinitely many simple closed geodesics on its surface. 
They can be naturally enumerated by pairs of coprime integers~$n> m> 1$
with two additional cases~$(1,0)$ and $(1,1)$. 
The $(n,m)$-geodesic is a broken line  with~$4(n+m)$ vertices, its length 
tends to infinity as~$m\to \infty$. Are there other polyhedra possessing 
this property?
The answer depends on convexity.  
We give an elementary proof that among convex polyhedra only
disphenoids admit arbitrarily long closed simple geodesics. 
For  non-convex polyhedra, this is not true. We present a counterexample with  the corresponding 
 polyhedron being a union of seven equal cubes. Several open problems are formulated.

\bigskip

\smallskip

\begin{flushright}

\end{flushright}

\end{abstract}

\begin{center}

\large{\textbf{1. Introduction}}	
\end{center}
\bigskip

Geodesics on a surface are locally shortest curves.  
They play a role of straight lines on the plane, however,
may possess different properties. For example,   
geodesics can have many points of intersections or they can be closed curves. Closed geodesics have important applications in the dynamical systems and differential geometry. In 1905 A.\,Poincar\'e conjectured that 
on every convex smooth surface  in~$\re^3$ there exist at least three closed and simple (i.e., non self-intersecting) geodesics. In other words,   there are at least three ways to 
pull the rubber band  around a pebble so that it will not slip off. 
The conjecture was proved in 1930 by L.\,Lusternik and L.\,Schnirelman~\cite{LS} and later
\cite{B, F}  
it was shown that  each convex smooth surface actually admits infinitely many 
 closed  simple  geodesics. Why the hypothesis was formulated only for three geodesics? 
 The number ``three''  was, probably, 
 explained by the general belief that  the ellipsoid with different lengths of axes 
 admits exactly three geodesics, which are its intersections with  three planes of symmetry.  
 This mistake was, in particular, made by D.\,Hilbert and S.\,Cohn-Vossen in ``Geometry and the imagination''~\cite{HC}. Elementary examples 
 of the ``fourth'' geodesic on an ellipsoid were constructed much later, see~\cite{K}. 
 
 Geodesics on polyhedra are systematically studied since  1990s
 and have become a subject of an extensive literature by now. 
 The geometrical properties of geodesics are
 analysed for regular polyhedra in~\cite{AA, AAH, DDTY, DL, FF, F2} and for 
 general polyhedra in~\cite{G1, G2, P1, P2}. The relation between 
  geodesics and  billiards was addressed in~\cite{F1, F2}; the 
  geodesics in the spherical and Lobachevsky geometry are considered in~\cite{BS1, BS2}.

 Most of convex polyhedra do not admit closed simple geodesics at all. This fact does not 
 contradict the Lusternik-Schnirelman theorem, since the surface of a polyhedron is not smooth but only piecewise smooth.  
If such a geodesic exists, then it divides the surface   into 
two domains and  for each of them, the sum of face angles over the vertices in that 
domain must be a multiple of~$2\pi$.
This  is a direct consequence of the Gauss-Bonnet theorem, although there are more elementary proofs~\cite{G1, G2}. For example, a tetrahedron can have a simple closed geodesic only if the sum of angles at two of its vertices is~$2\pi$. That is why a right triangular pyramid distinct from a regular tetrahedron   has no simple closed geodesics at all. On the other hand, all regular polyhedra do possess closed simple geodesics~\cite{FF}. The geodesics on the cube 
 have only three different length, those on the regular octahedron have two. 
  Here are the squares of these lengths (the length of the
edge is one):
\smallskip 

{\footnotesize Cube: $\quad 16,\, 18,\, 20$;}

{\footnotesize Octahedron: $\quad 9,\, 12$;}

{\footnotesize Icosahedron: $\quad 25,\, 27,\, 28$;}

{\footnotesize Dodecahedron: $\quad 27 + 18\phi\, , \, 28 + 20\phi\, , \, 
 ; \, 29 + 18\phi\, ; \, 29 + 19\phi\, , \, 25 + 25\phi$, }
\smallskip 

\noindent where $\phi= \frac{\sqrt{5}+1}{2}$ is the golden ratio. 
A regular tetrahedron admits arbitrarily long  simple closed geodesics.  Say, a regular tetrahedron with the  edge 1 cm possesses a closed 
geodesic longer than 1 km. Without self-intersections!  The same is true for every 
disphenoid, an isosceles 
tetrahedron,  whose faces are equal acute-angled triangles~\cite{FF, P2}. 
A natural question arises 
if there are other polyhedra that admit arbitrarily long geodesics? 
For convex polyhedra, the answer turns out to be 
negative. We proved it in~\cite{P2}, where we also  conjectured 
that the disphenoid is unique in this property among not only polyhedra but all convex surfaces. 
The conjecture  was proved in 2018 by  A.\,Akopyan and A.\,Petrunin~\cite{AP}. What can be said about non-convex polyhedra? Can they have arbitrary long geodesics?

Since  geodesics on polyhedra are broken  lines with vertices on edges, 
 they can be studied by tools of elementary geometry. 
This is what we are going to do in this paper. 

First, we recall a compete classification of 
closed simple geodesics on the surface of disphenoid. This implies, in particular, that 
the disphenoid admits arbitrarily long geodesics. 
Then we show that no other convex polyhedra possess this property.
To this end, we present a new elementary proof.  
Finally we solve this problem for the non-convex case. 

Several open problems
will be formulated. One of them concerns geodesics in spherical and Lobachevsky spaces. Do they have polyhedra (possibly, non-convex) with infinitely many simple closed geodesics?  
It is known that regular tetrahedra  in those spaces do not possess this property~\cite{BS1, BS2}.  
\bigskip 

\newpage

\begin{center}
\large{\textbf{2. Geodesics on the surface of a polyhedron}}	
\end{center}
\bigskip 

We begin with necessary definitions. Let~$\Tor$ denote the unit circle. 
\smallskip 

{\em A closed geodesic} on a surface  $S$ is a 
rectifiable curve 
 $\, \ell: \, \Tor \to  S\, $ that possesses the following property: 
for every  $t \in \Tor$, there exists  $\varepsilon >0$ such that for every $\tau_1, \tau_2 \in (t - \varepsilon, t + \varepsilon)$, 
the arc of $\ell$ is the shortest way on~$S$ between the points
$\ell(\tau_1)$ and $\ell(\tau_2)$. We always assume that a geodesic is simple, i.e., 
without self-intersections.

If~$P$ is a polyhedron, then every geodesic is a broken line
with vertices on edges. Indeed, on each face, 
 $\ell$ must be  a straight line as a unique locally-shortest 
curve on the plane. We call the vertices of~$\ell$ {\em nodes}, not to confuse them with 
vertices of~$P$.  The sides of~$\ell$ are {\em links}; several  consecutive 
links form a {\em chain}. 

Every node~$V$ is either a vertex of~$P$ or an interior point of an edge 
with the ``reflection property'': two links passing through~$V$
make equal angles with that edge (all angles are oriented, the links are directed 
in succession). This means that after unfolding the two faces onto the plane these links become 
collinear. 

\begin{figure}[h!]
\includegraphics{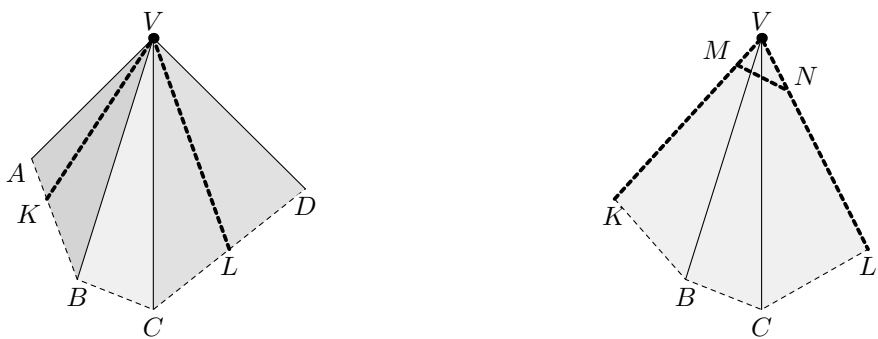}
 	\caption{{\footnotesize A geodesic on a polyhedron cannot pass through the vertices}}
 	\label{f.10}
 \end{figure}

Can a geodesic~$\ell$ pass through a vertex~$V$ of the polyhedron~$P$?  
If the  sum of angles of faces at~$V$  is less than $2\pi$, 
then it cannot. 
Otherwise, if~$\ell$ contains~$V$, then we cut the surface of the polyhedral angle~$V$
into two parts along the links~$KV, LV$ that meet at~$V$. The sum of angles in one of those parts is less than~$\pi$. We unfold this part onto the plane and take two close points on different 
links (points $M$ and $N$ on fig.~\ref{f.10}). 
Then the shortest way between those points is the segment~$MN$, which does not pass through~$V$. 

Thus, a geodesic does not contain vertices at which the  sum of face angles  
 is less than $2\pi$. In particular,  {\em if the polyhedron is convex, then 
 the geodesic cannot pass though the vertices.}

For the convex polyhedron, the aforementioned properties are sufficient
to define a geodesic: 
\begin{prop}\label{p.10}
 A simple closed curve on a convex polyhedron is geodesic if and only if the 
following hold: 

1) this is a broken line with nodes on edges; 

2) it has the reflection property at every node;  

3) it does not pass through the vertices.
\end{prop}

This fact enables us to study geodesics on polyhedra by elementary geometrical tools. 
The convexity is significant: for non-convex polyhedra, 
a simple closed geodesic 
can contain vertices and can go along edges. 
\smallskip 

\begin{figure}[h!]
\includegraphics{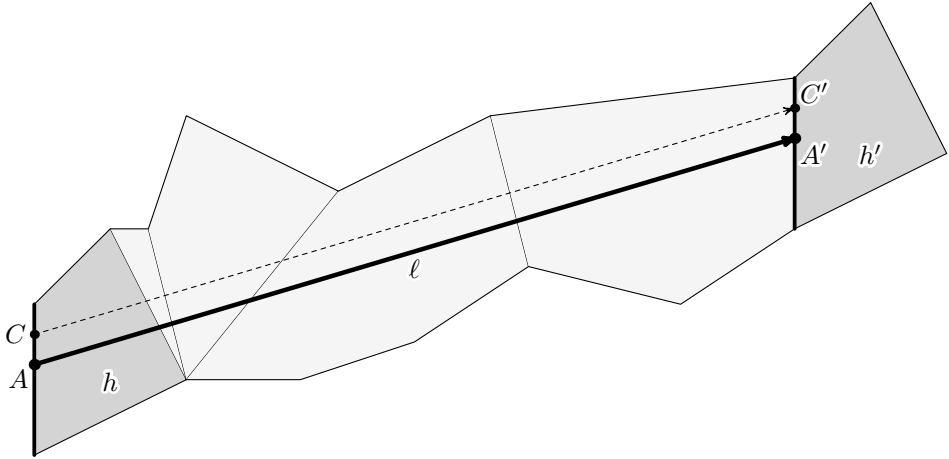}
 	\caption{{\footnotesize Equivalent geodesics~$AA'$ and~$CC'$ on the diagram}}
 	\label{f.60}
 \end{figure}

If one consequently  develops the faces of~$P$ onto the plane along the geodesic~$\ell$, 
then~$\ell$  becomes a straight line. 
This line will be denoted by the same symbol~$\ell$. 
The unfolding starts at some face~$h$ and closes in
its  copy~$h'$. The polygon~$h'$ is a translation of~$h$ by a vector 
parallel to~$\ell$ of length equal to the perimeter of~$\ell$, fig.~\ref{f.60}
We call two geodesics {\em equivalent} if they intersect the same 
edges in the same order.

\begin{prop}\label{p.20}
 If a convex polyhedron admits a simple closed  
geodesic, then it has infinitely many equivalent geodesics of the same length.  
\end{prop}
{\tt Proof.} We unfold  the faces of the polyhedron onto the plane of the net along the geodesic~$\ell$. 
We begin at  some node~$A \in \ell$ and move along the geodesic unfolding the faces 
which~$\ell$ meets. Having passed the whole geodesic we arrive at the same node~$A$, which becomes a point~$A'$ on the net. Thus, the geodesic~$\ell$ becomes a segment~$AA'$  which does not pass through
the images of vertices on the net. Slightly move the point~$A$ along the edge and obtain 
a point~$C$ (fig.~\ref{f.60}). Then the segment~$CC'$  
connects two images of the same node and does not pass through the vertices, provided that~$C$ is close enough to~$A$.  Since the segments~$AC$ and~$A'C'$ are equal and parallel, 
$ACC'A'$ is a parallelogram. Therefore, $CC'$ is equal and parallel to~$AA'$. 
Thus, the segment~$CC'$ also represents a geodesic, which is equivalent to~$\ell$
and has the same length.

{\ \hfill $\Box$}
\smallskip


\bigskip 

\newpage

\begin{center}
\large{\textbf{3. Isosceles tetrahedron: arbitrarily long geodesics}}	
\end{center}
\bigskip 

The disphenoid possesses an exceptional property: it has the {\em full net} that 
fills the whole plane. 

\begin{figure}[h!]
\includegraphics{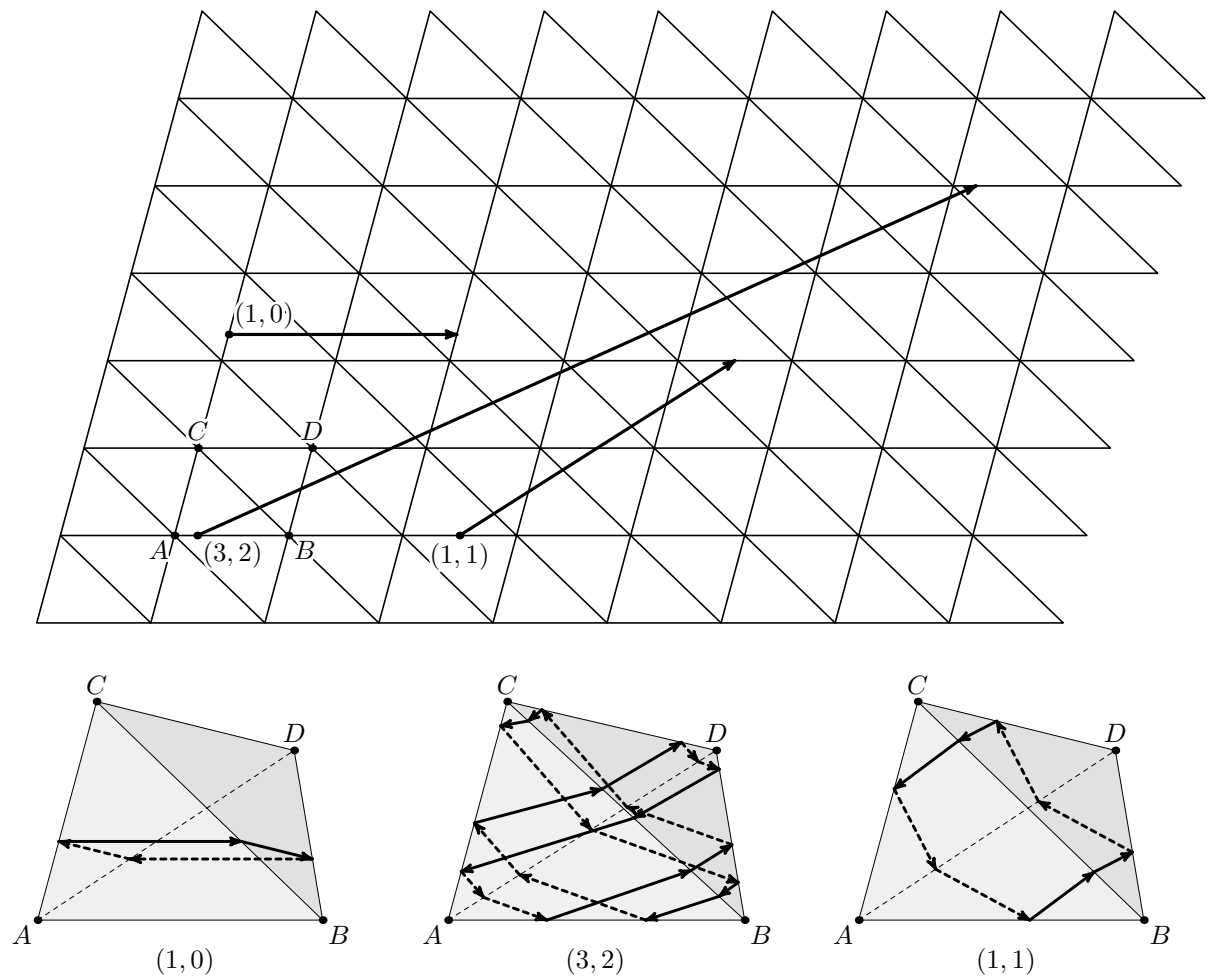}
 	\caption{{\footnotesize Geodesics of types $(1,0), \, (3,2)$
and~$(1,1)$} on the disphenoid}
\label{f.20}
 \end{figure}

Take a disphenoid~$ABCD$ and consider a triangular lattice 
generated by the face~$ABC$. This is the  full net  with knots corresponding to the vertices of~$ABCD$. This means that  one can label the vertices of
the lattice by the letters $A, B, C, D$ so that  for any development the
labelling of vertices of the tetrahedron will match the labelling of vertices of the
lattice,~see fig.~\ref{f.20}. Nothing like that can be done for the cube (or the
octahedron, or the icosahedron).

Introduce the system of coordinates on the plane with the origin at~$A$
and with the basis vectors~$2 AB, \, 2 AC$. Thus, 
$\, B\, =\, \bigl(\frac12, 0 \bigr), \, C\, =\, \bigl(0, \frac12
\bigr)$. Each point~$X$ on the surface of the tetrahedron has the full preimage on the plane: 
 $\, \bigl\{ \, \pm \, X \, + \, \bu\, ,\ \bu \in
\, \z^2\, \bigr\}$, where $\z^2$ is the set of integer vectors. 
In particular, the vertex $A$ corresponds to all integer points on the plane, 
$B \mapsto \bigl(n + \frac12 , m \bigr), \, C \mapsto \bigl(n  , m + \frac12\bigr), \,
D \mapsto \bigl(n + \frac12 , m + \frac12\bigr), \  m,n \in \z$.

Now we apply the full net to analyse geodesics on a disphenoid. 
Two geodesics on a polyhedron are called {\em isomorphic} if one of them 
becomes equivalent to the other after an automorphism 
(isometric map to itself) of the polyhedron. 
 
\begin{theorem}\label{th.30}
Consider the full net of the disphenoid.  Every straight line on this net  parallel to 
an integer vector and not containing the knots represents 
a closed simple geodesic. All  geodesics parallel to one vector are isomorphic 
and have the same length.
\end{theorem}
{\tt Proof.} Let~$\ell$ be parallel to 
an integer vector $\bu=(n,m)$. It may be assumed that  $n \ge m \ge 0$
and $n,m$ are co-prime, unless~$m=0$.
If the points~$X, Y \in \ell$ represent the same point on the tetrahedron, then they 
are either symmetric about a knot or $XY$ is an integer vector. 
The former is impossible, otherwise~$\ell$ contains a   knot (the center of symmetry);   
in the latter case the vector $XY$ must be an integer multiple of~$\bu$, since 
$\bu$ has co-prime coordinates. After the division, it may be assumed that~$XY = \bu$. 
Hence, the segment~$XY$ is folded onto a closed geodesic, which is simple 
since all interior points of that segment  represent different points of the surface of tetrahedron.

It remains to prove that parallel lines on the net produce  isomorphic geodesics. 
Among all straight lines on the net parallel to~$\ell$ and passing through knots,  
choose the closest ones from each side, call them~$\ell_1, \ell_2$. 
On each of those lines we take a segment which connects two knots and does not contain other knots. 
 Since $\ell_2$ is obtained from~$\ell_1$ by an integer translation 
  leaving the net invariant, 
 those two segments are equal. They are opposite sides of a parallelogram which 
 does not  contain 
 other knots since neither do the whole strip between~$\ell_1$ and~$\ell_2$
 (we denote this strip by~$\ell_1\ell_2$). The  vertices of this parallelogram present 
 four different vertices of the tetrahedron, otherwise, the midpoint of the 
 segment between the two corresponding knots would be a knot. 
 So, this parallelogram can be denoted as~$A'B'C'D'$, where 
 $X'$ represents a point~$X$ on the net.  
Any other line~$a$ parallel to~$\ell$ on the net can also be put into a strip
$a_1a_2$ whose sides contain vertices of a parallelogram~$A''B''C''D''$. 
 A suitable integer translate or 
symmetry about a knot maps the strip $\, \ell_1\ell_2\, $ to $\, a_1a_2$,
and the knot  $A'\, $ to $\, A''\, $. Consequently, we can assume  that 
 $\, a\, $ lies in the same strip  $\,
\ell_1\ell_2$. Since there are no knots inside this strip, 
it follows that the lines~$\ell$ and $a$ intersect the same 
edges in the same order, i.e., these geodesics  are isomorphic. 

{\ \hfill $\Box$}
\smallskip

\begin{remark}\label{r.30}{\em In the proof of Theorem~\ref{th.30} we actually show that  
every geodesic on a surface of a dispheniod is simple. 
Indeed, within each face a geodesic makes a family of 
parallel lines.
}
\end{remark}

The correspondence between closed simple geodesics on the isosceles tertrahedron 
and the straight lines  with rational directions in the full triangular net is
 well-known, see~\cite{CF, FF}. The isomorphism of all trajectories parallel to  a given rational vector was observed in~\cite{P1}. This gives a combinatorial 
classification of simple closed geodesics on disphenoid.
Each of them is unfolded to a straight line  parallel to an integer vector~$(n,m)$ 
with co-prime components or with~$n=1, m=0$. Vice versa, every such a vector produces a class of isomorphic geodesics of equal lengths~$|\ell| = \sqrt{n^2 + m^2 + mn}$.  
It is shown easily that~$\ell$  
has  $\, m, n\, $ and $\, n+m\, $ nodes on the edges  $\, AB, AC\, $
 and $\, BC$ respectively, i.e.,~$4(n+m)$ nodes in total.

Figure~\ref{f.20} presents the geodesics of types $(n,m) \, = \, (1,0), \, (3,2)$
and~$(1,1)$. The geodesic~$(3,2)$ has $4(3+2) = 20$ links.

As we have mentioned, a right triangular pyramid has no closed simple geodesics, provided 
that the lateral  edge is not equal to the side of the base. Hence, if we slightly move 
the vertex of the regular tetrahedron up along the altitude, then all its 
(infinitely many!) geodesics will slip off.  

\begin{remark}\label{r.40}
{\em Actually, our classification  holds for geodesics on an {\em arbitrary}  
tetrahedron, not only on an isosceles one. Namely, if we associate 
vertices of an arbitrary tetrahedron~$\Delta$ to the vertices 
of a disphenoid, then {\em each closed simple geodesic 
on~$\Delta$ is isomorphic to one on the disphenoid}~\cite{P1}.
This means that if~$\Delta$ admits a simple closed geodesic, then
it intersects the edges of the tetrahedron in the same order 
as one of the geodesics on the disphenoid does. 
Therefore, each simple closed geodesic on the tetrahedron has a combinatorial type $(n,m)$ described above. 
 Every non-isosceles tetrahedron has only finitely many 
(maybe zero) non-equivalent geodesics, we prove this is the next section. 
All of them have the aforementioned combinatorial types. 

Let us recall that a generic tetrahedron does not possess closed simple geodesics at all~\cite{G2}. On the other hand, it
is still not known (although looks plausible) that a generic tetrahedron possesses at least one
closed geodesic,  simple or not.}
\end{remark}

\bigskip

\begin{center}
\large{\textbf{4. Other convex polyhedra}	}
\end{center}
\bigskip 

Do there exist convex polydedra apart from disphenoids 
that admit infinitely many non-isomorphic 
closed simple geodesics or, equivalently, arbitrary long geodesics? 
At first sight, such examples can be obtained merely by 
slight perturbations of a disphenoid. For instance,  by cutting off  its vertices
(replacing them by small triangular faces) one obtains a polyhedron which is very close to 
the dispenoid. Clearly, it may inherit  a lot of non-equivalent geodesics from the disphenoid, provided that the cuts are small enough. However, it may not have infinitely many ones. It turns out that actually it never does. 
The following theorem established in~\cite{P2}  asserts that the existence of arbitrary long closed simple geodesics 
is an exceptional property of disphenoids. Here we give an elementary proof.

\begin{theorem}\label{th.10}
If a convex polyhedron has arbitrarily long simple closed geodesics, then it
is  a disphenoid. 
\end{theorem}

The roadmap of the proof is the following. 
Every simple closed geodesic splits the surface into two pieces~$\Gamma$ and~$\tilde \Gamma$. 
We show that if~$\ell$ is sufficiently long, then  unfolding~$Gamma$ to the plane  we meet a very long and narrow straight strip 
with almost parallel lateral sides. Then we prove that this strip ends by hooking a vertex 
of~$P$.  
This implies that~$\Gamma$ is folded from one  strip and contains two vertices. 
The same is true for~$\tilde \Gamma$. Hence, 
$P$ has four vertices, i.e., is a tetrahedron.   Then we show that this tetrahedron is
isosceles.   
The proof is split into four steps. 
\medskip 
   
\noindent \textbf{1.} {\em The domain~$\Gamma$ possesses a net,  a part of which is a  
straight strip with many nodes on the sides.}
\smallskip

Every edge of the polyhedron~$P$ can contain several nodes of 
the geodesic~$\ell$. At most two of them are {\em extreme}, i.e., 
closest to 
the ends of the edge, all others  are {\em intermediate}. Since the total number of 
extreme nodes does not exceed the double number of edges of~$P$, there are pieces 
of the geodesic~$\ell$ that contain only intermediate nodes~$A_1, \ldots , A_N$. 
The number~$N$ can be as large as we want,  
provided that~$\ell$ 
is long enough. Every node~$A_i$ lies on an edge of~$P$, where 
it has two neighbours. 
For one of them~(call it~$B_i$), the segment~$A_iB_i$ 
 lies in~$\Gamma$. Thus, we obtain a sequence of nodes~$B_1, \ldots , B_N$
 such that each segment~$A_iB_i$ is a part of one edge of~$P$. 
 We call the link~$A_iA_{i+1}$ {\em incorrect } if the 
 nodes~$B_i$ and~$B_{i+1}$ are not connected by a link. 
 In this case they are ends of two different links going  on the  same  face 
 on one half-plane with respect to~$A_iA_{i+1}$. If~$\ell$ contains one incorrect link 
 connecting two nodes on the same  edges, then it lies on the other 
 half of that face and hence, there are no other incorrect links. Therefore,  each pair of edges of~$P$ produces at most two incorrect links.
 Thus, the total number of incorrect links of~$\ell$ is bounded above by the double number of pairs of edges of~$P$.  
 Consequently, there are arbitrarily long sequences of the nodes~$A_i$ without 
incorrect links. We may assume that all links in the chain~$A_1 \cdots A_N$
are correct.  Hence,~$B_1 \cdots  B_N$ is also a chain and 
the piece of~$\Gamma$ between those chains is a belt winded around~$P$. 
Denote this belt by~$S_N$. 
Unfolding it onto the plane   
 obtain a straight strip on the net bounded by the segments~$A_1A_N$ and~$B_1B_N$. 
 \medskip 
 
 \noindent \textbf{2.} {\em Those stripes are long and narrow, their lateral sides are almost parallel.} 
\smallskip 

Thus, for arbitrary large~$N$, there is a belt~$S_N$ 
bounded by two chains $A_1\cdots A_N$ and $B_1\cdots B_N$ of correct links.  
We denote by~$M$ the maximal number of edges of~$P$ starting at one vertex and by~$d$ the smallest distance 
between disjoint edges of~$P$. Let us show that the perimeter of 
every subchain~$A_{i} \cdots A_{i+M}$ of~$M$ links
is at least~$d$.  This will imply that the chain~$A_1\cdots A_N$ has a big length, 
provided that~$N$ is large enough.   
Suppose  the contrary: the perimeter of some subchain
~$A_{i} \cdots A_{i+M}$ is smaller than~$d$.  
Hence, every two nodes~$A_j, A_k$ of this subchain belong 
to neighbouring edges of~$P$. Indeed, if they are from disjoint edges, then 
the length~$A_jA_k$ is not less than the distance between those edges, 
which is at least~$d$. This is impossible since the 
arc of~$\ell$ between~$A_j$ and~$A_k$ is shorter than~$d$. 
Therefore, all edges containing the  nodes~$A_{i},  \ldots , A_{i+M}$  
are neighbouring, i.e., start from the same vertex of~$P$. 
Denote this vertex by~$V$. Let $m$ be the total number of edges starting at~$V$.
Since~$m \le M$, we see that two of those nodes must be  on the same edge. 
Let them be~$A_i$ and $A_{i+k}, \, k\le m$, and all nodes between them are on different edges. This means that 
 the line~$\ell$ 
walks around~$k$ edges  intersecting them at the nodes~$A_{i}, \ldots , A_{i+k-1}$
and  comes back to the first edge (containing~$A_i$) at the node~$A_{i+k}$.
This implies that~$\ell$ passes through all edges starting at~$V$, i.e, $\, k=m$.

 If the segment $VA_{i+m}$ is shorter than~$VA_{i}$, then 
the  line~$\ell$ must go further from~$A_{i+m}$ 
around the same edges  and intersect them 
at nodes located closer to the vertex~$V$, otherwise~$\ell$ has self-intersections. 
This way $\ell$ will make infinitely many loops and never ends. If $VA_{i+m}$ is longer 
 than~$VA_{i}$, then we make the same conclusion moving along~$\ell$ in the opposite direction:  from the node~$A_i$ to $A_{i-1}$, etc.

  Thus, taking~$N$ large enough one can make the belt~$S_N$
arbitrarily long. Hence, the quadrangle~$A_1A_N B_N B_1$
 on the plane of the net has long sides~$A_1A_{N}$ and  $B_1B_{N}$. 
 Its area, however,  does not exceed the surface area of~$P$. 
   Hence, those sides are almost parallel 
and close to each other. This means that the angle and the distance between them
are both very small.  
 
  \medskip 
 
 \noindent \textbf{3.} {\em 
 The domain~$\Gamma$ contains exactly two vertices, each has the 
 sum of  face angles close to~$\pi$.} 
\smallskip 

\begin{figure}[h!]
\centerline{\includegraphics{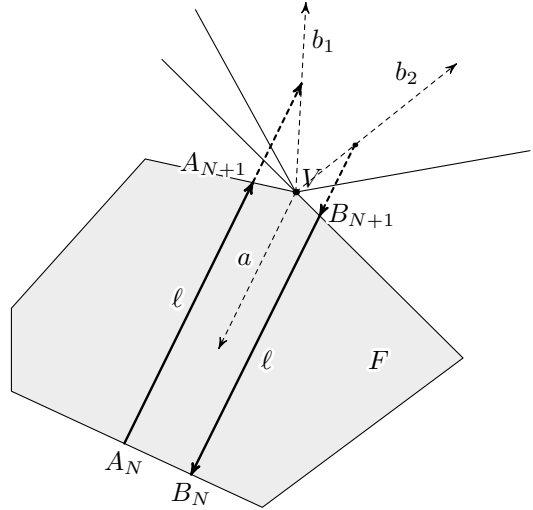}}
 	\caption{{\footnotesize The geodesic~$\ell$ hooks the vertex~$V$}}
\label{f.30}
 \end{figure}

The domain~$\Gamma$ contains a belt~$S_N$
 with arbitrarily long and ``almost parallel'' lateral sides. 
This belt ends when two links~$A_{N}A_{N+1}$ and $B_{N}B_{N+1}$ arrive at 
the nodes~$A_{N+1}$ and~$B_{N+1}$ located on different edges
of a single face~$F$ of~$P$.  The
part of~$F$  between those links is in~$\Gamma$, hence, 
it does not contain other nodes.   
The segment~$A_{N+1}B_{N+1}$ is very short, otherwise~$S_N$ will have a large area, provided 
that it is long enough. Therefore, 
the arc~$A_{N+1}B_{N+1}$ of the boundary of~$F$  contains a unique vertex~$V$. 
We draw a ray~$a$ on~$F$ from this vertex  parallel to~$A_{N}A_{N+1}$. 
Then draw a ray~$b$ from~$V$  on another face of~$P$ such that 
$a$ and $b$ split the surface of the solid angle~$V$ into two parts
with equal sum of face angles of faces. Each sum is smaller than~$\pi$, denote it by~$\pi -
\varepsilon$.  
Making a  cut along~$b$ (after which~$b$ becomes two rays~$b_1$ and~$b_2$) and unfolding the faces around~$V$ onto the plane 
we obtain~$\angle \,a b_1\, = \, \angle \,a  b_2 \, = \, \pi - \varepsilon$. Since the links~$A_{N}A_{N+1}$ and~$B_{N}B_{N+1}$  are 
almost parallel to~$a$ (the first one is exactly parallel), it follows 
that the straight lines containing those links on the net intersect the rays~$b_1$ and $b_2$ 
respectively (fig.~\ref{f.30}).

  This means that the arc of the geodesic~$\ell$ connecting the points~$A_{N+1}$ and $B_{N+1}$  intersects all the edges starting at~$V$  and intersects the ray~$b$. 
  Thus, the line~$\ell$ makes a loop around~$V$
changing the direction from the vector~$A_{N}A_{N+1}$ to the vector~$B_{N+1}B_{N}$. 
Those vectors are almost parallel, consequently, the sum of face angles at the vertex~$V$ is close to~$\pi$. 

Thus, the line~$\ell$ ``hooks up'' the vertex~$V$. Extending the  belt~$S_N$ 
in the opposite direction we prove that it hooks up another vertex and thus the domain~$\Gamma$
closes. Therefore,~$\Gamma$ contains exactly two vertices, 
the  sum of  face angles at each of them is  close to~$\pi$. 

 \bigskip 
 
  \noindent \textbf{4.} {\em The polyhedron~$P$ is a disphenoid.} 
\smallskip 
 
 Since we can apply the same line of reasoning to both domains~$\Gamma$ and $\tilde \Gamma$, we see that 
 each of them contains two vertices.  Hence,~$P$ is a tetrahedron. 
 The sum of  face angles at each vertex can be made arbitrarily close to~$\pi$, hence, it is equal to~$\pi$. It remains to invoke the following well-known fact: 
 {\em a tetrahedron is isosceles if and only if the sum of face angles at each vertex is equal to~$\pi$}. This completes the proof. 
\bigskip

{\ \hfill $\Box$}
\smallskip

\medskip

\begin{center}
\large{\textbf{5. Non-convex polyhedra}}	
\end{center}

\medskip

Can a non-convex polyhedron~$P$ have infinitely many simple closed geodesics? The proof of Theorem~\ref{th.10} is violated because of the only reason: 
there may be vertices, at which the  sum of face angles exceeds~$2\pi$. If~$P$ does not possess such vertices, then all face angles of its faces are less than $\pi$ and hence, 
all faces are convex polygons. Thus, all faces of~$P$ are convex and the  sum of face angles 
at each vertex is less than~$2\pi$. In this case the 
 proof of Theorem~\ref{th.10} is directly 
extended to~$P$. Thus, we obtain:

\begin{prop}\label{p.40}
If a non-convex polyhedron has the sum of face angles at each vertex 
less  than $2\pi$, then it does not admit arbitrarily long closed simple geodesics. 
\end{prop}

If the polyhedron has a vertex~$V$ with the sum of face angles at least~$2\pi$, then
the geodesic can pass through~$V$ and, moreover, can go  from~$V$ along an edge. 
Even if the geodesic forms a narrow strip, it now may not 
hook the vertex~$V$.  We are going to see that this is enough to make Theorem~\ref{th.10}
 false for non-convex polyhedra. 
 
 Before we formulate the main result, we clarify 
 one aspect. Let a surface~$S$ bound some compact domain~$G \subset \re^3$. 
  If a geodesic~$\ell$ on~$S$ is locally shortest not only 
 among curves laying on~$S$ but also among all curves not intersecting the 
 interior of~$G$, then~$\ell$ is called {\em strong}. In other words, 
 every short arc of a strong geodesic is the shortest path 
 between its ends among all paths not intersecting~${\rm int}\, G$. 
 This stronger property possesses the physical  interpretation of a geodesic as 
 a robber band on the surface. We call the geodesics possessing this property  
 {\em strong}. Non-strong geodesics may hot have a form of an elastic band. For convex surfaces,  all geodesics are strong.

Proposition~\ref{p.40} is true for all geodesics, not necessarily strong.  
 Now we present a non-convex polyhedron with 
 arbitrarily long {\em strong} geodesics. The construction consists of seven equal cubes. 
 
 
 \begin{theorem}\label{th.40}
 There exists a non-convex polyhedron which possesses arbitrary long 
closed geodesics that are simple and strong. 
 \end{theorem}
 \enlargethispage{2\baselineskip} 
 
  \noindent \textbf{ The main idea of the construction.} Assemble four equal cubes face-to-face 
 as in fig.~\ref{f.40} (left). 
 Choose an arbitrary point~$M$ on the edge~$EF$ and connect it 
 with the vertex~$D$ by a 
 geodesic line~$M K D$.  Then draw a segment $AL$  parallel
 to~$MK$ with~$L\in EF$.  We make two key observations: 
\smallskip

\begin{figure}[h!]
\includegraphics{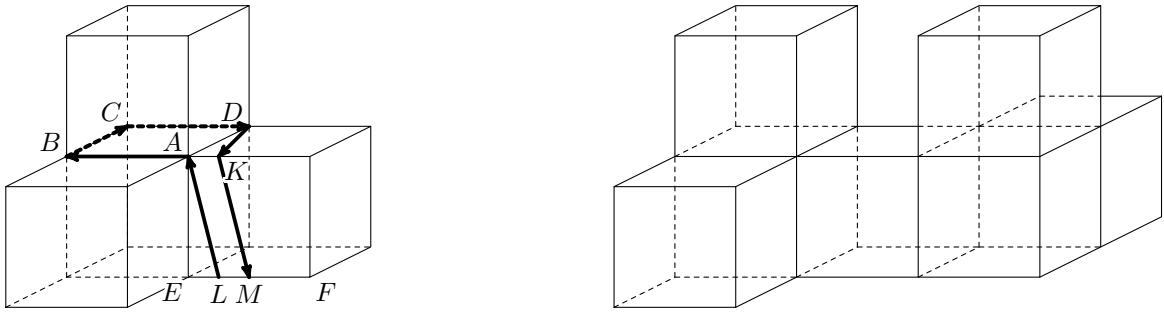}
 	\caption{{\footnotesize Construction of the polyhedron}}
\label{f.40}
 \end{figure}  

\smallskip 

 \noindent \textbf{1)} {\em The broken line $LABCDKM$ is a geodesic connecting the points~$M$ and $L$}. This is shown in a straightforward manner by unfolding.

 \noindent \textbf{2)} {\em The strip between~$KM$ and~$AL$ can be made arbitrary 
 narrow} by moving the point~$M$ closer to~$E$.  
  \smallskip

\begin{figure}[h!]
\includegraphics{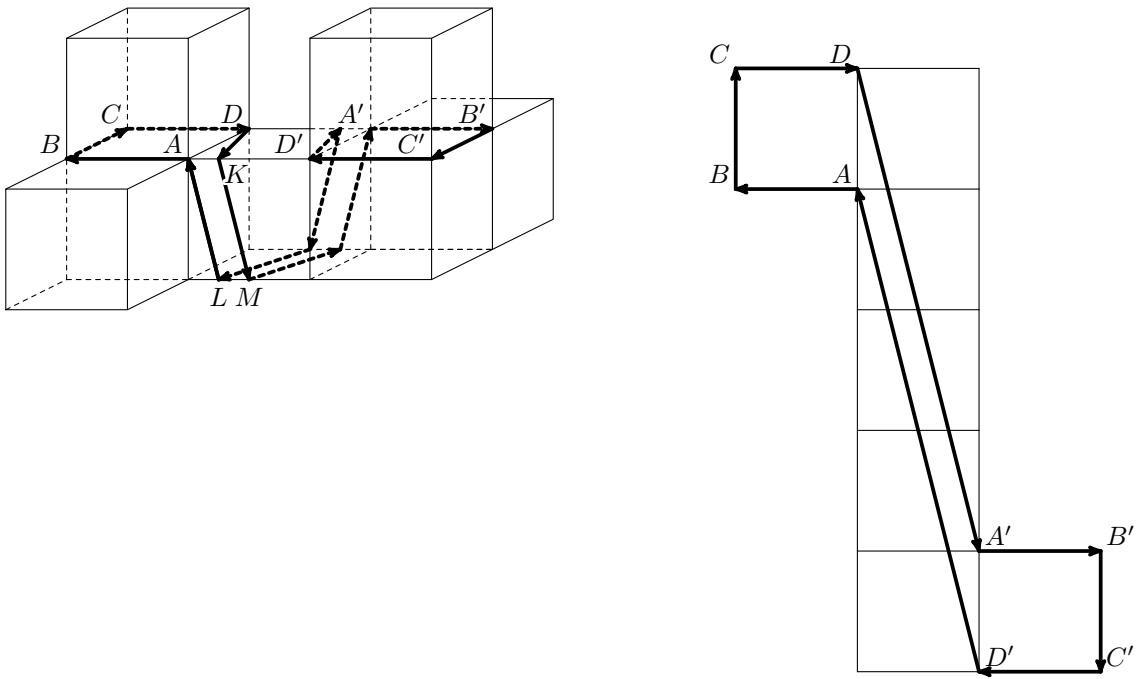}
 	\caption{{\footnotesize Arbitrarily long geodesics}}
\label{f.50}
 \end{figure}

    \noindent \textbf{ The  construction.} We put two figures from~fig.~\ref{f.40} (left) together so that 
    they share a common middle cube. Thus we obtain the polyhedron~$P$ that consists of seven cubes, fig.~\ref{f.40} (right). 
    The line~$\ell$ makes a loop around the upper cube along the 
    path~$ABCD$, then makes several turns  around the middle cube and then 
    makes  a symmetric loop~$A'B'C'D'$. Since~$\ell$
    can form an arbitrarily narrow strip, it can make arbitrarily many turns 
    around the middle cube without self-intersections.     
    Now we present the formal construction (fig.\ref{f.50}).

     Consider the net obtained by~$k$
    consecutive unfolding  of the middle cube onto the plane. It consists of  $4k+1$ square faces,
    the first and the last ones represent the same face~$ADA'D'$.  
    Two  parallel segments $AD'$ and $DA'$ and two letters~$\Pi$ 
    on the upper and lower squares form the  
    image of the closed geodesic~$\ell$. 
    Taking~$k$ large enough we make~$\ell$ as long as needed. 
     
\smallskip

\bigskip

\begin{center}
\large{\textbf{6. Discussion and open problems}}	
\end{center}
\bigskip

The  geodesics  constructed in the proof of Theorem~\ref{th.40} contain  
whole edges of the polyhedron. Can this property be avoided? 
\smallskip

  \noindent \textbf{Problem~1}. {\em Can a non-convex polyhedron admit arbitrarily long 
  closed simple geodesic not containing edges?} 
\smallskip 

Similarly to the convex case, one can show that each of the two parts of the surface bounded by a geodesic possesses long and narrow straight strips with almost parallel sides. 
The only difference is  that now the line~$\ell$ may hook several vertices. 
 If the sum of face angles at those vertices is  $\pi k$, where $k \ge 3$ is an odd number, 
then the sides of the strip will be parallel. The question is can those sides be arbitrarily close to each other, which is necessary for long geodesics?  

\smallskip

  \noindent \textbf{Problem~2}. {\em What is the minimal number of vertices of a non-convex polyhedron that admits arbitrarily long 
  closed simple geodesics ?} 
\smallskip 

 \noindent  For example, can it have 5, 6, 7 vertices?

\smallskip

It is known that regular tetrahedra  in spherical and Lobachevsky spaces   have only finitely many closed simple geodesics~\cite{BS1, BS2}. A natural question arises are there
polyhedra (may be non-convex) in those spaces with infinitely many geodesics? 
 \smallskip 
 
  \noindent \textbf{Problem~3}. {\em  Do there exist  polyhedra (possibly, non-convex) 
  in spherical and Lobachevsky spaces that admit  infinitely many simple closed geodesics?}

 \end{document}